\documentclass{article}
\usepackage{amssymb}
\usepackage{amsfonts}

\input{tcilatex}
\begin{document}

\begin{center}
\textbf{The} \textbf{\bigskip local and global} \textbf{dynamics of a cancer
tumor growth and chemotherapy treatment model }

\textbf{Veli Shakhmurov}

Department of Mechanical Engineering, Okan University, Akfirat, Tuzla 34959
Istanbul, Turkey,

E-mail: veli.sahmurov@okan.edu.tr\ 

\textbf{A. Maharramov}

Department of Biophysics Yeditepe University

\bigskip E-mail: amaharramov@yeditepe.edu.tr

\textbf{Bunyad Shahmurzada}

Kocaeli University, Faculty of Medicine, Umuttepe Campus, Izmit, Kocaeli,
Turkey, E-mail: shahmurzadab@yahoo.com

\textbf{Abstract}
\end{center}

In this paper, we studied phase-space analysis of a certain mathematical
model of tumor growth with an immune responses and chemotherapy therapy.
Mathematical modelling of this process is viewed as a potentially powerful
tool in the development of improved treatment regimens. Mathematical
analysis of the model equations with multipoint initial condition, regarding
nature of equilibria, local and global stability have been investigated. We
studied some features of behavior of one of three-dimensional tumor growth
models with dynamics described in terms of densities of three cells
populations: tumor cells, healthy host cells and effector immune cells. We
found sufficient conditions, under which trajectories from the positive
domain of feasible multipoint initial conditions tend to one of equilibrium
points. The addition of a drug term to the system can move the solution
trajectory into a desirable basin of attraction. We show that the solutions
of the model with a time-varying drug term approach can be evaluated more
fruitful way and down to earth style from the point of practical importance
than the solutions of the system without drug treatment, in the condition of
stimulated immune processes, only.


\begin{center}
\textbf{1. Introduction}
\end{center}

Beginning with this article, we intend to attempt to investigate the
problems of mathematical and biological approaches to the modeling of cancer
growth dynamics processes and operations, enlisting a new parameter as
chemotherapeutic agent to the system, consisted of host, tumor and immune
cells. The mathematical processing of the model is based on nonlinear
property of cancer growth, when it concerns the foundation of the model's
logistic part. This approach appears very convenient in description of
unexpected dynamics in the processes of growth in response to changing
reactions of the system to different concentrations of immune cells and drug
application at the different stages of cancer growth development $\left[ 1-13%
\right] $. Taking into account all the complex processes, nonlinear
mathematical models can be estimated capable of compensation and
minimization the inconsistencies between different mathematical models
related to cancer growth-anticancer factor affections. The elaboration of
mathematical non-spatial models of the cancer tumor growth in the broad
framework of tumor immune interactions studies is one of intensively
developing areas in the modern mathematical biology, see works $[1-7]$. Of
course, for the development of a powerful mathematical model on cancer
immunotherapy, it is required first of all, an understanding of the
mechanisms governing the dynamics of tumor growth.\ One of the main reasons
for creation of non-spatial dynamical models of a multiphase nature is
related to the fact that they are described by a system of ordinary
differential equations, which can be efficiently investigated by powerful
methods of qualitative theory of ordinary differential equations and
dynamical systems theory. In this paper, we examine the dynamics of a cancer
growth model with drug interaction proposed in $[5],$ taking account the
system as a multiphase structure, i.e. in dynamics: 
\[
\dot{T}=r_{1}T\left( 1-k_{1}^{-1}T\right) -a_{12}NT-a_{13}TI-g_{1}\left(
u\right) T, 
\]%
\begin{equation}
\dot{N}=r_{2}N\left( 1-k_{2}^{-1}N\right) -a_{21}NT-g_{2}\left( u\right) N,%
\text{ }  \tag{1.1}
\end{equation}%
\[
\dot{I}=s+\frac{r_{3}IT}{k_{3}+T}-a_{31}IT-d_{3}I-g_{3}\left( u\right) I, 
\]%
\[
\dot{u}+d_{2}u\left( t\right) =\upsilon \left( t\right) 
\]%
with multipoint initial condition 
\begin{equation}
T\left( t_{0}\right) =T_{0}+\dsum\limits_{j=1}^{m}\alpha _{1j}T\left(
t_{j}\right) \text{, }N\left( t_{0}\right)
=N_{0}+\dsum\limits_{j=1}^{m}\alpha _{2j}N\left( t_{j}\right) \text{, } 
\tag{1.2}
\end{equation}%
\[
I\left( t_{0}\right) =I_{0}+\dsum\limits_{j=1}^{m}\alpha _{3j}I\left(
t_{j}\right) \text{, }t_{0}\in \left[ 0,\right. \left. \delta \right) \text{%
, }t_{j}\in \left( 0,\delta \right) ,\text{ }t_{j}>t_{0},\text{ }u\left(
t_{0}\right) =u_{0}, 
\]%
where $T=T\left( t\right) ,$ $N=N\left( t\right) $, $I=I\left( t\right) $
denote the densities of tumor cells, healthy host cells and the effector
immune cells respectively, at the moment $t,$ $k_{i}>0,$ $\alpha _{ij}$ are
real numbers, $m$ is a natural number such that 
\[
T\left( t_{0}\right) >0,\text{ }N\left( t_{0}\right) >0,\text{ }I\left(
t_{0}\right) >0 
\]%
and 
\begin{equation}
g_{i}\left( u\right) \geq 0\text{, }g_{i}\left( 0\right) =0,\text{ }%
\lim\limits_{u\rightarrow \infty }\text{ }g_{i}\left( u\right) =\text{\ }%
a_{i}>0\text{, }i=1,2,3.\text{ }  \tag{1.3}
\end{equation}

For

\begin{equation}
g_{i}\left( u\right) =a_{i}\left( 1-e^{-\nu _{i}u}\right) ,\text{ }\nu
_{i}>0,\text{ }i=1,2,3  \tag{1.4}
\end{equation}%
we generalize the case that has been derived in $\left[ 5\right] .$

The source of the immune cells is considered to be outside of the system so
it is reasonable to assume a constant influx rate $s$, furthermore, in the
absence of any tumor, the cells will die off at a per capita rate $d_{3}$,
resulting in a long-term population size of $s/d_{3}$ cells, $u(t)$ denotes
the amount of drug at the tumor site at time $t$, this is determined by the
dose given $\upsilon (t)$, and a per capita decay rate of the drug once it
is injected, here it is assumed that the drug kills all types of cells, but
that the kill rate differs for each type of cell, with the response curve in
all cases given by%
\[
g(u)=\left( g_{1}(u),g_{2}(u),g_{3}(u)\right) . 
\]%
For case of $\left( 1.4\right) $, $g(u)$ is the fractional cell kill for a
given amount of drug $u$, at the tumor site, this decay rate incorporates
all pathways of elimination of the drug, by $a_{1}$, $a_{2}$ and $a_{3}$
denoted the three different response coefficients, here%
\[
a=\left( a_{1},a_{2},a_{3}\right) ,\text{ }\nu \text{ }=\left( \nu _{1},\nu
_{2},\nu _{3}\right) . 
\]%
The first term of the first equation corresponds to the logistic growth of
tumor cells, in the absence of any effect from other cells populations, with
the growth rate of $r_{1}$ and maximum carrying capacity $k_{1}$. The
competition between host cells and tumor cells $T\left( t\right) $ which
results in the loss of the tumor cells population is given by the term $%
a_{12}NT$. Next, the parameter $a_{13}$ refers to the tumor cell killing
rate by the immune cells $I\left( t\right) $. In the second equation, the
healthy tissue cells also grow logistically with the growth rate of $r_{2}$
and maximum carrying capacity $k_{2}$. We assume that the cancer cells
proliferate faster than the healthy cells which gives $r_{1}>r_{2}$. The
tumor cells also inactivate the healthy cells at the rate of $a_{21}$. The
third equation of the model describes the change in the immune cells
population with time $t.$ The first term of the third equation illustrates
the stimulation of the immune system by the tumor cells with tumor specific
antigens. The rate of recognition of the tumor cells by the immune system
depends on the antigenicity of the tumor cells. The model of the recognition
process is given by the rational function which depends on the number of
tumor cells with positive constants $r_{3}$ and $k_{3}$. The immune cells
are inactivated by the tumor cells at the rate of $a_{31}$ as well as they
die naturally at the rate $d_{3},$ here we suppose that the constant influx
of the activated effector cells into the tumor microenvironment is zero. We
suppose that the constant influx $s$ of the activated effector cells into
the tumour microenvironment is zero.

Therein, note that, the nonlinear dynamic systems and references studied
e.g. in $[14-15]$. One of main aims is derivation of sufficient conditions
under which the possible biologically feasible dynamics is local and
globally stable, and converges to one of equilibrium points. Since these
equilibrium points have a biological sense, we notice that understanding
limit properties of dynamics of cells populations based on solving problems $%
(1.1)-\left( 1.2\right) $ may be of an essential interest for the prediction
of health conditions of a patient without a treatment, when the data (e.g.
the status of blood cells shown above) that determines the condition of the
patient are compared at various times $t_{0},t_{1},...,t_{m}$ and
correlated. Note that the local and global stability properties of $\left(
1.1\right) $ with the classical initial condition were studied in $\left[ 8%
\right] $ and $\left[ 9\right] $, respectively. We prove that all orbits are
bounded and must converge to one of several possible equilibrium points.
Therefore, the long-term behavior of an orbit is classified according to the
basin of multipoint attraction in which it starts. By scaling $%
x_{1}=Tk_{1}^{-1}$, $x_{2}=Nk_{2}^{-1}$, $x_{3}=Ik_{3}^{-1}$, $\tilde{t}%
=r_{1}t$ in $\left( 1.1\right) -\left( 1.2\right) $ and omitting the tilde
notation we obtain the multipoint initial value problem (IVP) 
\[
\dot{x}_{1}=x_{1}\left( 1-x_{1}\right)
-a_{12}x_{1}x_{2}-a_{13}x_{1}x_{3}-g_{1}\left( u\right) x_{1}, 
\]%
\begin{equation}
\dot{x}_{2}=r_{2}x_{2}\left( 1-x_{2}\right) -a_{21}x_{1}x_{2}-g_{2}\left(
u\right) x_{2},  \tag{1.5}
\end{equation}%
\[
\dot{x}_{3}=\frac{r_{3}x_{1}x_{3}}{x_{1}+k_{3}}%
-a_{31}x_{1}x_{3}-d_{3}x_{3}-g_{3}\left( u\right) x_{3},\text{ }t\in \left[
0,\right. \left. T\right) , 
\]

\[
\dot{u}+d_{2}u\left( t\right) =\upsilon \left( t\right) ,\text{ }u\left(
t_{0}\right) =u_{0}, 
\]

\begin{equation}
x_{1}\left( t_{0}\right) =x_{10}+\dsum\limits_{j=1}^{m}\alpha
_{1j}x_{1}\left( t_{j}\right) \text{, }x_{2}\left( t_{0}\right)
=x_{20}+\dsum\limits_{j=1}^{m}\alpha _{2j}x_{2}\left( t_{j}\right) \text{, }
\tag{1.6}
\end{equation}%
\[
x_{3}\left( t_{0}\right) =x_{30}+\dsum\limits_{j=1}^{m}\alpha
_{3j}x_{3}\left( t_{j}\right) \text{, }t_{0}\in \left[ 0,\right. \left.
\delta \right) \text{, }t_{j}\in \left( 0,T\right) ,\text{ }t_{j}>t_{0}, 
\]%
where $\alpha _{ij}$ are real numbers and $m$ is a natural number such that 
\begin{equation}
x_{j0}+\dsum\limits_{k=1}^{m}\alpha _{jk}x_{j}\left( t_{k}\right) \geq 0%
\text{, }j=1,2,3.  \tag{1.7}
\end{equation}%
By solving the problem 
\begin{equation}
\dot{u}\left( t\right) +d_{2}u\left( t\right) =\upsilon \left( t\right) 
\text{, }u\left( t_{0}\right) =u_{0}  \tag{1.8}
\end{equation}%
we get%
\begin{equation}
u\left( t\right) =e^{t_{0}}u_{0}+\dint\limits_{t_{0}}^{t}e^{-\left( t-\tau
\right) }\upsilon \left( \tau \right) d\tau ,\text{ }t\in \left( 0,\left. T%
\right] ,\right.  \tag{1.9}
\end{equation}%
i.e. the sustem $\left( 1.5\right) $ is equivalent to the following 
\[
\dot{x}_{1}=x_{1}\left( 1-x_{1}\right)
-a_{12}x_{1}x_{2}-a_{13}x_{1}x_{3}-g_{1}\left( \upsilon \right) x_{1}, 
\]%
\begin{equation}
\dot{x}_{2}=r_{2}x_{2}\left( 1-x_{2}\right) -a_{21}x_{1}x_{2}-g_{2}\left(
\upsilon \right) x_{2},  \tag{1.10}
\end{equation}%
\[
\dot{x}_{3}=\frac{r_{3}x_{1}x_{3}}{x_{1}+k_{3}}%
-a_{31}x_{1}x_{3}-d_{3}x_{3}-g_{3}\left( \upsilon \right) x_{3},\text{ }t\in %
\left[ 0,\right. \left. T\right) , 
\]%
where $g_{i}\left( \upsilon \right) =f_{i}\left( u\right) ,$ here $u$ is
defined by $\left( 1.9\right) .$ Note that, for $\alpha _{j1}=\alpha
_{j2}=...\alpha _{jm}=0$ the problem $\left( 1.5\right) -\left( 1.6\right) $
turns to be the Cauchy problem%
\[
\dot{x}_{1}=x_{1}\left( 1-x_{1}\right)
-a_{12}x_{1}x_{2}-a_{13}x_{1}x_{3}-g_{1}\left( u\right) x_{1}, 
\]%
\begin{equation}
\dot{x}_{2}=r_{2}x_{2}\left( 1-x_{2}\right) -a_{21}x_{1}x_{2}-g_{2}\left(
u\right) x_{2},  \tag{1.11}
\end{equation}%
\[
\dot{x}_{3}=\frac{r_{3}x_{1}x_{3}}{x_{1}+k_{3}}%
-a_{31}x_{1}x_{3}-d_{3}x_{3}-g_{3}\left( u\right) x_{3},\text{ } 
\]%
\[
\dot{u}+d_{2}u\left( t\right) =\upsilon \left( t\right) ,\text{ }t\in \left[
0,T\right] , 
\]%
\[
x_{1}\left( t_{0}\right) =x_{10}\text{, }x_{2}\left( t_{0}\right) =x_{20},%
\text{ }x_{3}\left( t_{0}\right) =x_{30},\text{ }u\left( t_{0}\right) =u_{0},%
\text{ }t_{0}\in \left[ 0,\right. \left. T\right) . 
\]

\begin{center}
\textbf{2}.\textbf{\ Notations and background.}
\end{center}

Consider the multipoint IVP for nonlinear equation%
\begin{equation}
\frac{du}{dt}=f\left( u\right) ,\text{ }t\in \left[ 0,T\right] ,  \tag{2.1}
\end{equation}%
\[
u\left( t_{0}\right) =u_{0}+\dsum\limits_{k=1}^{m}\alpha _{k}u\left(
t_{k}\right) \text{, }t_{0}\in \left[ 0,\right. \left. T\right) \text{, }%
t_{k}\in \left( 0,T\right) ,\text{ }t_{k}>t_{0} 
\]%
in a Banach space $X$, where $\alpha _{k}$ are complex numbers, $m$ is a
natural number and $u=u\left( t\right) $ is a $X-$valued function. Note
that, for $\alpha _{1}=\alpha _{2}=...\alpha _{m}=0$ the problem $\left(
2.1\right) $ becomes the following local Cauchy problem%
\begin{equation}
\frac{du}{dt}=f\left( u\right) ,\text{ }u\left( t_{0}\right) =u_{0},\text{ }%
t\in \left[ 0,T\right] ,\text{ }t_{0}\in \left[ 0,\right. \left. T\right) 
\text{.}  \tag{2.2}
\end{equation}

For $u_{0}\in X$ let $\bar{B}_{r}\left( u_{0}\right) $ denotes a closed ball
in $X$ with radius $r$ centered at $u_{0}$, i.e., 
\[
\bar{B}_{r}\left( u_{0}\right) =\left\{ u\in X:\left\Vert u-u_{0}\right\Vert
_{X}\leq r\right\} . 
\]

From $\left[ 19\right] $ we have

\textbf{Theorem 2.1. }Let $X$ be a Banach space. Suppose $f:X\rightarrow X$
satisfies local Lipschitz condition on $\bar{B}_{r}(\upsilon _{0})\subset $ $%
X$, i.e.%
\[
\left\Vert f\left( u\right) -f\left( \upsilon \right) \right\Vert _{X}\leq
L\left\Vert u-\upsilon \right\Vert _{X} 
\]%
for each $u$, $\upsilon \in \bar{B}_{r}(\upsilon _{0})$ and there exists $%
\delta >0$ such that 
\[
t_{k}\in O_{\delta }\left( t_{0}\right) =\left\{ t\in \mathbb{R}:\left\vert
t-t_{0}\right\vert <\delta \right\} , 
\]%
where 
\[
\upsilon _{0}=u_{0}+\dsum\limits_{k=1}^{m}\alpha _{k}u\left( t_{k}\right) . 
\]

Moreover, let 
\[
M=\sup\limits_{u\in \bar{B}_{r}(\upsilon _{0})}\left\Vert f\left( u\right)
\right\Vert _{X}<\infty . 
\]

Then,\ problem $\left( 2.1\right) $ has a unique continuously differentiable
local solution $u(t)$, for $t\in O_{\delta }\left( t_{0}\right) $, where $%
\delta \leq \frac{r}{M}.$

\textbf{Theorem 2.2. }Let $X$ be a Banach space. Suppose that $%
f:X\rightarrow X$ satisfies global Lipschitz condition, i.e.%
\[
\left\Vert f\left( u\right) -f\left( \upsilon \right) \right\Vert _{X}\leq
L\left\Vert u-\upsilon \right\Vert _{X} 
\]%
for each $u$, $\upsilon \in X$. Moreover, let 
\[
M=\sup\limits_{u\in X}\left\Vert f\left( u\right) \right\Vert _{X}<\infty . 
\]

Then\ problem $\left( 2.1\right) $ has a unique continuously differentiable
local solution $u(t)$, for $\left\vert t-t_{0}\right\vert <\delta $, where $%
\delta \leq \frac{r}{M}.$

Let $X$ be a Banach space. $w\in X$ is called a critical point (or
equilibria point) for the equation $\left( 2.1\right) $ if $f\left( w\right)
=0.$

We denote the solution of the problem $\left( 2.1\right) $ by%
\[
\phi \left( t,u_{0}\right) =\phi \left( t,u\left( t_{0}\right) ,u\left(
t_{1}\right) ,...,u\left( t_{m}\right) \right) . 
\]

\textbf{Definition 2.1. }Let $u_{0}\in X$. \textbf{\ }A critical point $w\in
X$ of the equation $\left( 2.1\right) $ is called a positive multiphase
attractor if there exists a neighbourhood $O_{w}\subset X$ of $w$ such that
the relation 
\[
u_{0}=u\left( t_{0}\right) -\dsum\limits_{k=1}^{m}\alpha _{k}u\left(
t_{k}\right) \subset O_{w}\text{ for }t_{0}\in \left[ 0,\right. \left.
T\right) \text{, }t_{k}\in \left( 0,T\right) ,\text{ }t_{k}>t_{0} 
\]%
implies $\lim\limits_{t\rightarrow \infty }u\left( t\right) =w.$

\textbf{Definition 2.2. }Assume $w\in X$ is a multiphase attractor of $%
\left( 2.1\right) .$ A set of $u_{0}\in X$ with a property that for solution 
$\phi \left( t,u_{0}\right) $ of $\left( 2.1\right) $ we have $%
\lim\limits_{t\rightarrow \infty }u\left( t\right) =w,$ is called a domain
of multiphase attractor (domain of multiphase asymptotic stability, or
multiphase basin) of $w.$\ 

\begin{center}
\textbf{3. The equilibria points, existence and local stability}
\end{center}

The equilibria points of the system $(1.3)$ are obtained by solving the
system of isocline equations%
\[
x_{1}\left( 1-x_{1}\right) -a_{12}x_{1}x_{2}-a_{13}x_{1}x_{3}-g_{1}\left(
\upsilon \right) x_{1}=0, 
\]%
\begin{equation}
r_{2}x_{2}\left( 1-x_{2}\right) -a_{21}x_{1}x_{2}-g_{2}\left( \upsilon
\right) x_{2}=0,  \tag{3.1}
\end{equation}%
\[
\frac{r_{3}x_{1}x_{3}}{x_{1}+k_{3}}-a_{31}x_{1}x_{3}-d_{3}x_{3}-g_{3}\left(
\upsilon \right) x_{3}=0. 
\]%
Since we are interested in biologically relevant solutions of $\left(
3.1\right) ,$ we find sufficient conditions under which this system have
positive solutions.

It is clear\textbf{\ }to see that\textbf{\ }the points $E_{0}\left(
0,0,0\right) $, $E_{1}\left( \gamma ,0,0\right) $, $E_{2}\left( 0,\delta
,0\right) $ are equilibria points for the system of $\left( 1.11\right) $,
where%
\begin{equation}
\gamma =1-g_{1}\left( \upsilon \right) ,\text{ }\delta =1-\frac{g_{2}\left(
\upsilon \right) }{r_{2}}.  \tag{3.2}
\end{equation}

\textbf{Remark 3.1.} It is clear that\textbf{\ }the points $E_{0}\left(
0,0,0\right) $, $E_{1}\left( \gamma ,0,0\right) $, $E_{2}\left( 0,\delta
,0\right) $ are biologically feasible equilibria points for the system $%
\left( 1.10\right) $.

Let%
\[
R_{+}^{3}=\left\{ x\in R^{3}\text{: }x_{i}>0\text{, }i=1,2,3\right\} . 
\]

\textbf{Remark 3.2.} (1) consider the equilibrium points $E_{0}\left(
0,0,0\right) ;$ for $E_{0}$ three type cell populations are zero; (2) for
points $E_{1}\left( \gamma ,0,0\right) $ tumor cells have survived but
normal and immune cells are zero, this case can be called as "dead" case;
(3) $E_{2}\left( 0,\delta ,0\right) $-tumor-free and immune free case; in
this category, normal cells have survived but tumor and immune cells are
zero.

We now, will derive that the linearized matrices of the system $\left(
1.11\right) $ for equilibria points $E_{0}\left( 0,0,0\right) $, $%
E_{1}\left( \gamma ,0,0\right) ,$ $E_{2}\left( 0,\delta ,0\right) $ are
following:%
\begin{equation}
A_{0}=\left[ 
\begin{array}{ccc}
\gamma & 0 & 0 \\ 
0 & d_{22} & 0 \\ 
0 & 0 & c_{33}%
\end{array}%
\right] ,\text{ }A_{1}=\left[ 
\begin{array}{ccc}
d_{11} & -a_{12}\gamma & -a_{13}\gamma \\ 
0 & d_{22} & 0 \\ 
0 & 0 & d_{33}%
\end{array}%
\right] ,\text{ }A_{2}=\left[ 
\begin{array}{ccc}
c_{11} & 0 & 0 \\ 
-a_{21}\delta & c_{22} & 0 \\ 
0 & 0 & c_{33}%
\end{array}%
\right] ,  \tag{3.3}
\end{equation}%
where 
\[
d_{11}=1-2\gamma -g_{1}\left( \upsilon \right) \text{, }d_{22}=r_{2}-a_{21}%
\gamma -g_{2}\left( \upsilon \right) ,\text{ } 
\]%
\begin{equation}
\text{ }d_{33}=\left( \frac{r_{3}}{\gamma +k_{3}}-a_{31}\right) \gamma
-d_{3}-g_{3}\left( \upsilon \right) ,\text{ }c_{11}=1-a_{12}\delta
-g_{1}\left( \upsilon \right) \text{,}  \tag{3.4}
\end{equation}%
\[
c_{22}=r_{2}\left( 1-2\delta \right) -a_{21}\delta -g_{2}\left( \upsilon
\right) \text{, }c_{33}=-d_{3}-g_{3}\left( \upsilon \right) . 
\]

\begin{center}
\textbf{4. local stability analysis of equilibria points }
\end{center}

In this section we show the following \ result:

\textbf{Theorem 4.1. (}1) $E_{0}\left( 0,0,0\right) $ is a locally
asymptotically stable point if $g_{1}\left( \upsilon \right) >1,$ $%
g_{2}\left( \upsilon \right) >r_{2}$ and $E_{0}$ is an unstable point if $%
g_{1}\left( \upsilon \right) <1,$ $g_{2}\left( \upsilon \right) <r_{2};$ (2) 
$E_{1}\left( \gamma ,0,0\right) $ is a locally asymptotically stable point
for the linearized system of $\left( 1.11\right) $ when $d_{ii}<0$ for $%
i=1,2,3,$ it is an unstable point when $d_{ii}>0$; (3) $E_{2}\left( 0,\delta
,0\right) $ is a locally asymptotically stable point for linearized system
of $\left( 1.11\right) $ when $c_{ii}<0$ for $i=1,2,$ it is an unstable
point if $c_{ii}>0$, where $\gamma $, $\delta $ were defined by $\left(
3.2\right) $ and $d_{ii}$, $c_{ii}$, $i=1,2,3$ defined by $\left( 3.4\right)
.$

\textbf{Proof. }Indeed, the eigenvalues of the matrix $A_{0}$ are $\gamma ,$ 
$d_{22},$ $c_{33};$ eigenvalues of the matrix $A_{1}$ are $d_{ii}$ and
eigenvalues of the matrix $A_{2}$ are $c_{ii}.$ Hence, by $\left[ \text{5,
Theorem 8.12}\right] $ we obtain the assertions.

\textbf{Theorem 4.2. }Let $g_{1}\left( \upsilon \right) >1$ and $g_{2}\left(
\upsilon \right) >r_{2}$. Then, the dimension of the stable manifold $%
W_{0}^{+}$ and unstable manifold $W_{0}^{--}$are given, respectively, by 
\begin{equation}
\text{Dim }W_{0}^{+}\left( E_{0}\left( 0,0,0\right) \right) =1,\text{ Dim }%
\left( W_{0}^{--}E_{0}\left( 0,0,0\right) \right) =1.  \tag{4.1}
\end{equation}

\textbf{Proof. }Let we solve the the following matrix equation%
\begin{equation}
A_{0}x=\lambda x,  \tag{4.2}
\end{equation}%
i.e. consider the system of homogenous linear equation 
\[
\left( 1-g_{1}\left( \upsilon \right) -\lambda \right) x_{1}=0,\text{ }%
\left( r_{2}-g_{2}\upsilon -\lambda \right) x_{2}=0_{,}\text{ }-\left(
d_{3}+g_{3}\left( \upsilon \right) +\lambda \right) x_{3}=0. 
\]%
By solving of $\left( 4.2\right) $ we get that eigenspaces corresponding to
eigenvalues $\lambda _{1}=\gamma ,$ $\lambda _{2}=d_{22},$ $\lambda
_{3}=\theta $ are respectively, the follwing: 
\[
B_{01}=\left\{ x\in R^{3}:x=\left( a,0,0\right) \right\} , 
\]%
\[
B_{02}=\left\{ x\in R^{3}:x=\left( 0,a,0\right) \right\} ,\text{ }%
B_{03}=\left\{ x\in R^{3}:x=\left( a,b,0\right) \right\} 
\]%
where $a$ is any real number, i.e. we obtain $\left( 4.1\right) .$

In a similar way we obtain

\textbf{Theorem 4.3. }Let $g_{1}\left( \upsilon \right) <1$ and $g_{2}\left(
\upsilon \right) <r_{2}$. Then, the dimension of the hyperbolic saddle
manifold $W^{0}$ is given by 
\[
\text{Dim }W^{0}\left( E_{0}\left( 0,0,0\right) \right) =1. 
\]

\textbf{Theorem 4.4. }Let $d_{ii}<0$, $i=1,2,3$. Then, the dimension of the
stable manifold $W_{1}^{+}$ and unstable manifold $W_{1}^{--}$are given,
respectively, by 
\begin{equation}
\text{Dim }W_{0}^{+}\left( E_{1}\left( \gamma ,0,0\right) \right) =1,\text{
Dim }W_{0}^{--}\left( E_{0}\left( \gamma ,0,0\right) \right) =1.  \tag{4.3}
\end{equation}

\textbf{Proof. }Let we solve the the following matrix equation%
\begin{equation}
A_{1}x=\lambda x,  \tag{4.4}
\end{equation}%
i.e. consider the system of homogenous linear equation 
\[
\left( d_{11}-\lambda \right) x_{1}-a_{12}\gamma x_{2}=0,\text{ }\left(
d_{22}-\lambda \right) x_{2}=0_{,}\text{ }\left( d_{33}-\lambda \right)
x_{3}=0. 
\]%
By solving of $\left( 4.4\right) $ we get that eigenspaces corresponding to
eigenvalues $\lambda _{1}=d_{11},$ $\lambda _{2}=d_{22},$ $\lambda
_{3}=d_{33}$ are respectively, the follwing: 
\[
B_{11}=\left\{ x\in R^{3}:x=\left( a,0,0\right) \right\} , 
\]%
\[
B_{12}=\left\{ x\in R^{3}:x=\left( b,a,0\right) \right\} ,\text{ }%
B_{13}=\left\{ x\in R^{3}:x=\left( 0,0,a\right) \right\} , 
\]%
where $a$ is any real number and 
\[
b=\frac{\left( 1-2\gamma -r_{2}+g_{2}\left( \upsilon \right) \right) }{%
a_{12}\gamma }, 
\]%
i.e. we obtain $\left( 4.3\right) .$

In a similar way we obtain

\textbf{Theorem 4.5. }Let $d_{ii}>0$. Then, the dimension of the hyperbolic
saddle manifold $W_{1}^{0}$ is given by 
\[
\text{Dim }W_{1}^{0}\left( E_{1}\left( \gamma ,0,0\right) \right) =1. 
\]

\textbf{Theorem 4.6. }Let $c_{ii}<0$, $i=1,2,3$. Then, the dimension of the
stable manifold $W_{2}^{+}$ and unstable manifold $W_{2}^{--}$are given,
respectively, by 
\begin{equation}
\text{Dim }W_{2}^{+}\left( E_{1}\left( 0,\delta ,0\right) \right) =1,\text{
Dim }W_{0}^{--}\left( E_{2}\left( 0,\delta ,0\right) \right) =1.  \tag{4.5}
\end{equation}

\textbf{Proof. }Let we solve the the following matrix equation%
\begin{equation}
A_{2}x=\lambda x,  \tag{4.6}
\end{equation}%
i.e. consider the system of homogenous linear equation 
\[
\left( c_{11}-\lambda \right) x_{1}=0,\text{ }-a_{21}\delta x_{1}+\left(
c_{22}-\lambda \right) x_{2}=0_{,}\text{ }\left( c_{33}-\lambda \right)
x_{3}=0, 
\]%
where $c_{ii}$ were defined by $\left( 3.4\right) .$ By solving of $\left(
4.6\right) $ we get that eigenspaces corresponding to eigenvalues $\lambda
_{1}=c_{11},$ $\lambda _{2}=c_{22},$ $\lambda _{3}=c_{33}$ are respectively,
the follwing: 
\[
B_{21}=\left\{ x\in R^{3}:x=\left( a,0,0\right) \right\} , 
\]%
\[
B_{22}=\left\{ x\in R^{3}:x=\left( 0,a,0\right) \right\} ,\text{ }%
B_{23}=\left\{ x\in R^{3}:x=\left( 0,0,a\right) \right\} , 
\]%
where $a$ is any real number. i.e. we obtain $\left( 4.5\right) .$

In a similar way we obtain

\textbf{Theorem 4.7. }Let $c_{ii}>0$. Then, the dimension of the hyperbolic
saddle manifold $W_{2}^{0}$ is given by 
\[
\text{Dim }W_{2}^{0}\left( E_{2}\left( 0,\delta ,0\right) \right) =1. 
\]

\bigskip \textbf{Definition 4.1.} A set $A\subset S$ is called a strong
multipoint attractor with respect to $S$ if%
\[
\limsup\limits_{t\rightarrow \infty }\rho \left( u\left( t\right) ,A\right)
=0, 
\]%
where $u\left( t\right) $ is an orbit such that $u\left( t_{0}\right)
-\dsum\limits_{k=1}^{m}\alpha _{k}u\left( t_{k}\right) \in $ $S$ and $\rho $
is the Euclidean distance function.

\textbf{Lemma 5.1.} $B$ is a strong multipoint attractor with respect to $%
R_{+}^{3}.$

\textbf{Proof. }The proof is done using standard comparison as in Theorem
3.1.

\begin{center}
\bigskip \textbf{5. Global stability of equilibria points }
\end{center}

\bigskip In this section, we will derive global stability condition of
equilibria points $E_{0}\left( 0,0,0\right) $, $E_{1}\left( \gamma
,0,0\right) $, $E_{2}\left( 0,\delta ,0\right) .$

Let 
\[
R_{+}^{3}=\left\{ x=\left( x_{1},x_{2},x_{3}\right) \in R^{n},\text{ }%
x_{k}\geq 0\right\} \text{, }\Omega _{K}=\left\{ x\in R^{3}\text{: }0\leq
x_{i}\leq K_{i}\text{, }i=1,2,3\right\} 
\]%
and 
\[
B_{r}\left( \bar{x}\right) =\left\{ x\in R^{3}\text{, }\left\Vert x-\bar{x}%
\right\Vert _{R^{3}}<r^{2}\right\} . 
\]

\textbf{Theorem 5.0. }Assume: (1) $g_{1}\left( \upsilon \right) >1$ and $%
g_{2}\left( \upsilon \right) >r_{2};$(2) $a_{31}k_{3}>r_{3}.$ Then the
system $\left( 1.11\right) $ is global asymptotically stabile at equilibria
point $E_{0}\left( 0,0,0\right) .$

\textbf{Proof. }Let $A_{0}$ be the linearized matrix with respect to
equilibria point $E_{0}\left( 0,0,0\right) ,$ i.e.%
\[
A_{0}=\left[ 
\begin{array}{ccc}
1-g_{1}\left( \upsilon \right) & 0 & 0 \\ 
0 & r_{2}-g_{2}\left( \upsilon \right) & 0 \\ 
0 & 0 & -d_{3}-g_{3}\left( \upsilon \right)%
\end{array}%
\right] .\text{ } 
\]

\bigskip We consider the Lyapunov equation 
\[
B_{0}A_{0}+A_{0}^{T}B_{0}=-I,\text{ }B_{0}=\left[ 
\begin{array}{ccc}
b_{11} & b_{12} & b_{13} \\ 
b_{21} & b_{22} & b_{23} \\ 
b_{31} & b_{32} & b_{33}%
\end{array}%
\right] . 
\]%
The equation above is reduced to linear system of algebraic equation with
respect to $b_{ij}.$ By solving this algebraic equation we get 
\[
b_{11}=\frac{1}{2\left( g_{1}\left( \upsilon \right) -1\right) },\text{ }%
b_{22}=\frac{1}{2\left( g_{2}\left( \upsilon \right) -r_{2}\right) },\text{ }%
b_{33}=\frac{1}{2\left( g_{3}\left( \upsilon \right) +d_{3}\right) }, 
\]%
\[
b_{ij}=0\text{, }i\neq j\text{, }i\text{, }j=1,2,3, 
\]%
i.e.,

\[
B_{0}=\left[ 
\begin{array}{ccc}
b_{11} & 0 & 0 \\ 
0 & b_{22} & 0 \\ 
0 & 0 & b_{33}%
\end{array}%
\right] . 
\]

Hence, 
\begin{equation}
P_{0}\left( \lambda \right) =\left\vert B_{0}-\lambda I\right\vert =\left(
b_{11}-\lambda \right) \left( b_{22}-\lambda \right) \left( b_{33}-\lambda
\right) =0.  \tag{5.0}
\end{equation}

By assumtion, eigenvalues $\lambda _{1}=b_{11},$ $\lambda _{2}=b_{22}$, $%
\lambda _{3}=b_{33}$ of the matrix $B_{0}$ are positive. So, the quadratic
function 
\begin{equation}
V_{0}\left( x\right)
=X^{T}B_{0}X=b_{11}x_{1}^{2}+b_{22}x_{2}^{2}+b_{33}x_{3}^{2}  \tag{$5.v_{0}$}
\end{equation}%
is a positive defined Lyapunov function candidate in the certain
neighborhood of $E_{0}\left( 0,0,0\right) .$ By $\left[ \text{12, Corollary
8.2}\right] $ we need now to determine a domain $\Omega _{0}$ about the
point $E_{1},$ where $\dot{V}_{0}\left( x\right) $ is negatively defined and
a constant $C$ such that $\Omega _{C}$ is a subset of $\Omega _{0}$. By
assuming $x_{k}\geq 0$, $k=1,2,3,$ we will find the solution set of the
following inequality 
\[
\dot{V}_{0}\left( x\right) =\dsum\limits_{k=1}^{3}\frac{\partial V_{0}}{%
\partial x_{k}}\frac{dx_{k}}{dt}=2b_{11}\left( 1-g_{1}\left( \upsilon
\right) \right) x_{1}^{2}+2b_{22}\left( r_{2}-g_{2}\left( \upsilon \right)
\right) x_{2}^{2}- 
\]%
\[
2b_{33}\left( d_{3}+g_{3}\left( \upsilon \right) \right)
x_{3}^{2}-2b_{11}x_{1}^{3}-2b_{22}x_{2}^{3}-2b_{11}g_{1}\left( \upsilon
\right) x_{1}^{2}- 
\]%
\[
2b_{11}x_{1}^{2}\left( a_{12}x_{2}+a_{13}x_{3}\right) -2b_{22}g_{2}\left(
\upsilon \right) x_{2}^{2}- 
\]%
\begin{equation}
2b_{22}a_{21}x_{2}^{2}x_{1}+2b_{33}x_{3}^{2}\left[ \frac{r_{3}x_{1}}{%
x_{1}+k_{3}}-a_{31}x_{1}\right] <0  \tag{5.1}
\end{equation}

For $x\in R_{+}^{3}$ we have 
\[
-2b_{11}x_{1}^{3}-2b_{22}x_{2}^{3}-2b_{11}x_{1}^{2}\left(
a_{12}x_{2}+a_{13}x_{3}\right) --2b_{22}a_{21}x_{2}^{2}x_{1}\leq 0. 
\]

Hence, in view of inequalities 
\begin{equation}
2ab\leq a^{2}+b^{2},\text{ }x_{1}^{2}+x_{2}^{2}\leq \left\Vert x\right\Vert
^{2},\text{ }x_{2}^{2}+x_{3}^{2}\leq \left\Vert x\right\Vert ^{2}  \tag{5.2}
\end{equation}%
for $x\in R_{+}^{3}$ we obtain that the inequality $\left( 5.1\right) $
holds if 
\[
2b_{11}\left( 1-g_{1}\left( \upsilon \right) \right) x_{1}^{2}+2b_{22}\left(
r_{2}-g_{2}\left( \upsilon \right) \right) x_{2}^{2}-2b_{33}\left(
d_{3}+g_{3}\left( \upsilon \right) \right) x_{3}^{2}\leq 0, 
\]%
\begin{equation}
\frac{r_{3}}{x_{1}+k_{3}}-a_{31}<0.  \tag{5.3 }
\end{equation}

By assumption (1) the first inequality of $\left( 5.3\right) $ are satisfied
for all $x\in R^{3}$ and the second \i nequlity holds by assumption (2). So, 
$\dot{V}_{0}\left( x\right) <0$ for $x\in R_{+}^{3},$ i.e. the point $E_{0}$
is global asymptotically stabile at equilibria point.

\bigskip \textbf{Theorem 5.1. }Assume: (1) $d_{11}<0,$ $d_{22}<0,$ $d_{33}<0;
$ (2); $-\left( d_{22}+d_{33}\right) <a_{21}\gamma ,$ 
\[
\gamma <-2\left( d_{11}+d_{22}\right) ,\text{ }-2d_{22}a_{12}^{2}\gamma
<a_{21}\left[ a_{12}^{2}\gamma +d_{11}\left( d_{11}+d_{22}\right) \right] ,
\]%
\[
\text{(3) }\frac{b_{22}}{2}\geq \frac{2b_{12}^{2}}{b_{11}},\text{ }\frac{%
b_{33}}{2}\geq \frac{2b_{13}^{2}}{b_{11}},\text{ }\frac{b_{33}}{2}\geq \frac{%
2b_{23}^{2}}{b_{22}};
\]

(4) 
\[
r_{2}>g_{2}\left( \upsilon \right) ,\text{ }\left( a_{12}+\gamma
d_{22}\right) a_{12}\gamma <d_{11}\left( d_{11}+d_{12}\right) . 
\]

Then the system $\left( 1.10\right) $ is global asymptotically stabile at
equilibria point $E_{1}\left( \gamma ,0,0\right) .$

\textbf{Proof. }Let $A_{1}$ be the linearized matrix with respect to
equilibria point $E_{1}\left( \gamma ,0,0\right) ,$ i.e. 
\[
A_{1}=\left[ 
\begin{array}{ccc}
d_{11} & -a_{12}\gamma  & -a_{13}\gamma  \\ 
0 & d_{22} & 0 \\ 
0 & 0 & d_{33}%
\end{array}%
\right] ,
\]%
where $d_{11}$ was defined by $\left( 3.4\right) .$

By assumption (1), $d_{ii}<0.$ We consider the Lyapunov equation 
\begin{equation}
B_{1}A_{1}+A_{1}^{T}B_{1}=-I,  \tag{5.4}
\end{equation}%
where 
\[
B_{1}=\left[ 
\begin{array}{ccc}
b_{11} & b_{12} & b_{13} \\ 
b_{21} & b_{22} & b_{23} \\ 
b_{31} & b_{32} & b_{33}%
\end{array}%
\right] . 
\]%
The equation $\left( 5.4\right) $ is reduced the linear system of algebraic
equation with respect to $b_{ij}$, by solving which we obtain 
\[
b_{11}=-\frac{1}{2d_{11}}\text{, }b_{12}=b_{21}=-\frac{a_{12}\gamma }{%
2d_{11}\left( d_{11}+d_{22}\right) },\text{ }b_{13}=b_{31}= 
\]%
\[
-\frac{a_{13}\gamma }{2d_{11}\left( d_{11}+d_{22}\right) },\text{ }b_{22}=%
\frac{1}{2d_{22}}\left( 2a_{12}b_{12}\gamma -1\right) , 
\]%
\begin{equation}
\text{ }b_{23}=\frac{\left( a_{13}b_{12}+a_{12}b_{13}\right) \gamma }{%
d_{22}+d_{33}},\text{ }b_{33}=-\frac{1}{2d_{13}}\left( 2a_{13}\gamma
b_{13}-1\right) .  \tag{5.5}
\end{equation}

Consider now, the quadratic function 
\[
V_{1}\left( x\right) =X^{T}B_{1}X=b_{11}\left( x_{1}-\gamma \right)
^{2}+b_{22}x_{2}^{2}+b_{33}x_{3}^{2}+ 
\]%
\[
2b_{12}\left( x_{1}-\gamma \right) x_{2}+2b_{13}\left( x_{1}-\gamma \right)
x_{3}+2b_{23}x_{2}x_{3}. 
\]%
It is clear to see that 
\[
V_{1}\left( x\right) =b_{11}\left( x_{1}-\gamma \right) ^{2}+2b_{12}\left(
x_{1}-\gamma \right) x_{2}+b_{22}x_{2}^{2}+2b_{13}x_{1}x_{3}+ 
\]%
\[
b_{33}x_{3}^{2}+2b_{23}x_{2}x_{3}=\frac{b_{11}}{2}\left( x_{1}-\gamma +\frac{%
2b_{12}}{b_{11}}x_{2}\right) ^{2}+\left[ \frac{b_{22}}{2}-\frac{2b_{12}^{2}}{%
b_{11}}\right] x_{2}^{2}+ 
\]%
\begin{equation}
\frac{b_{11}}{2}\left( x_{1}-\gamma +\frac{2b_{13}}{b_{11}}x_{3}\right) ^{2}+%
\left[ \frac{b_{33}}{2}-\frac{2b_{33}^{2}}{b_{11}}\right] x_{3}^{2}+ 
\tag{5.6}
\end{equation}%
\[
\frac{b_{22}}{2}\left( x_{2}+\frac{2b_{23}}{b_{22}}x_{3}\right) ^{2}+\left[ 
\frac{b_{33}}{2}-\frac{2b_{23}^{2}}{b_{22}}\right] x_{3}^{2}\geq 0, 
\]%
i.e. $V_{1}\left( x\right) $ is a positive defined Lyapunov function
candidate in neighborhood of $E_{1}\subset \Omega _{K}$ when the assumption
(3) hold. By $\left[ \text{12, Corollary 8.2}\right] $ we need now to
determine a domain $\Omega _{1}$ about the point $E_{1},$ where $\dot{V}%
_{1}\left( x\right) $ is negatively defined and a constant $C$ such that $%
\Omega _{C}$ is a subset of $\Omega _{1}$. By assuming $x_{k}\geq 0$, $%
k=1,2,3,$ we will find the solution set of the following inequality 
\begin{equation}
\dot{V}_{1}\left( x\right) =\dsum\limits_{k=1}^{3}\frac{\partial V_{1}}{%
\partial x_{k}}\frac{dx_{k}}{dt}=  \tag{5.7}
\end{equation}%
\[
2\left[ b_{11}\left( x_{1}-\gamma \right) +b_{12}x_{2}+b_{13}x_{3}\right]
x_{1}\left[ \left( 1-x_{1}\right) -a_{12}x_{2}-a_{13}x_{3}-g_{1}\left(
\upsilon \right) \right] + 
\]%
\[
2\left[ b_{12}\left( x_{1}-\gamma \right) +b_{22}x_{2}+b_{23}x_{3}\right]
x_{2}\left[ r_{2}\left( 1-x_{2}\right) -a_{21}x_{1}-g_{2}\left( \upsilon
\right) \right] + 
\]%
\[
2\left[ b_{23}x_{2}+b_{13}\left( x_{1}-\gamma \right) +b_{33}x_{3}\right]
x_{3}\left[ \frac{r_{3}x_{1}}{x_{1}+k_{3}}-a_{31}x_{1}-d_{3}-g_{3}\left(
\upsilon \right) \right] <0. 
\]

\bigskip By assumption (1) $b_{11},$ $b_{22}$, $b_{33}>0$, $b_{12}$, $%
b_{13}<0$, $b_{23}>0$. So, by assumption (2), some coefficients of terms $\
x_{1}x_{2},$ $x_{1}x_{2},$ $x_{3}^{3},\ x_{1}x_{2}^{2},$ $x_{1}^{2}x_{2},$ $%
x_{1}x_{2}x_{3}$ are negative. Hence, the estimate$\ \left( 5.7\right) $
holds if 
\begin{equation}
\eta _{11}x_{1}^{2}+\eta _{22}x_{2}^{2}++\eta _{33}x_{3}^{2}+2b_{11}\gamma 
\left[ a_{12}x_{1}x_{2}+a_{13}x_{1}x_{3}\right] <0,  \tag{5.8}
\end{equation}%
\[
\alpha _{1}x_{1}+\alpha _{2}x_{2}+\alpha _{3}x_{3}<0,
\]%
where, 
\[
\eta _{11}=2b_{11}\left[ 1+\gamma -g_{1}\left( \upsilon \right) \right] ,%
\text{ }\eta _{22}=2b_{12}\gamma +2r_{2}b_{22}-2b_{22}g_{2}\left( \upsilon
\right) ,
\]%
\[
\text{ }\eta _{33}=-2d_{3},\text{ }\alpha _{1}=2b_{11}\gamma \left[
g_{1}\left( \upsilon \right) -1\right] ,\text{ }\alpha _{2}=2b_{12}\gamma %
\left[ g_{2}\left( \upsilon \right) -r_{2}\right] ,
\]%
\[
\text{ }\alpha _{3}=2b_{13}\left[ g_{3}\left( \upsilon \right) +d_{3}\right]
.
\]%
\ By assumption (4) and in view of $\left( 5.5\right) ,$ $\eta _{ii}>0$ for $%
i=1,2,3.$ Moreover, by applying $\left( 5.2\right) $ we get that $\left(
5.8\right) $ hold if%
\begin{equation}
\mu _{1}\left( x_{1}-\gamma \right) ^{2}+\eta _{2}x_{2}^{2}+\mu
_{3}x_{3}^{2}<\mu _{1}\gamma ^{2},  \tag{5.9}
\end{equation}%
\[
\left( \alpha _{1}+2\mu _{1}\right) x_{1}+\alpha _{2}x_{2}+\alpha
_{3}x_{3}\leq 0,
\]%
where 
\[
\mu _{1}=\eta _{11}+b_{11}\gamma a_{12},\text{ }\mu _{2}=\eta
_{22}+b_{11}\gamma a_{13}.
\]%
Hence, $\dot{V}_{1}$ is negative defined on the domain $\Omega
_{1}=B_{r}\left( \bar{x}\right) \cap \Omega _{\gamma },$ if 
\begin{equation}
\left( x_{1}-\gamma \right) ^{2}+x_{2}^{2}+x_{3}^{2}<r^{2},  \tag{5.10}
\end{equation}%
\[
\left( \alpha _{1}+2\mu _{1}\right) x_{1}+\alpha _{2}x_{2}+\alpha
_{3}x_{3}\leq 0,
\]%
where%
\[
\bar{x}=\left( \gamma ,0,0\right) ,\text{ }r=\left( \frac{\mu _{0}}{\mu
_{1}\gamma ^{2}}\right) ^{\frac{1}{2}},\text{ }\mu _{0}=\max \left\{ \mu _{1}%
\text{, }\eta _{2}\text{, }\mu _{3}\text{ }\right\} ,\text{  }
\]%
\begin{equation}
r\leq \sqrt{K_{1}^{2}+K_{2}^{2}+K_{3}^{2}},\text{ }\Omega _{\gamma }=\left\{
x\in R_{+}^{3}\text{, }\left( \alpha _{1}+2\mu _{1}\right) x_{1}+\alpha
_{2}x_{2}+\alpha _{3}x_{3}\leq 0\right\} .  \tag{5.11}
\end{equation}%
i.e., the system $\left( 1.11\right) $ is global asymptotically\ stabile at $%
E_{1}\left( \gamma ,0,0\right) $.

\textbf{Remark 5.1. }In view of $\left( 5.5\right) $, the assumption (3) can
be realized as the condition on the coefficients of the system $\left(
1.11\right) .$

Let $c_{ii}$ be the numbers defined by $\left( 3.4\right) .$\ Now, we
consider the equilibria point $E_{2}\left( 0,\delta ,0\right) $ and prove
the following result

\textbf{Theorem 5.2. }Assume (1) $c_{11}<0,$ $\delta =1-\frac{g_{2}\left(
\upsilon \right) }{r_{2}}>0,$ $c_{22}<0$ (2) $g_{1}\left( \upsilon \right)
<1;$

(3) 
\[
g_{1}\left( \upsilon \right) +g_{2}\left( \upsilon \right)
+r_{2}>r_{2}\left( \delta -1\right) \text{, }g_{1}\left( \upsilon \right)
-1>\delta \left( a_{12}+1\right) . 
\]

Then the system $\left( 1.11\right) $ is global asymptotically stabile at
equilibria point $E_{1}\left( 0,\delta ,0\right) .$

\textbf{Proof. }Let $A_{2}$ be the linearized matrix with respect to
equilibria point $E_{2}\left( 0,\delta ,0\right) ,$ i.e.%
\[
A_{2}=\left[ 
\begin{array}{ccc}
c_{11} & 0 & 0 \\ 
-a_{21}\delta & c_{22} & 0 \\ 
0 & 0 & c_{33}%
\end{array}%
\right] . 
\]

Consider the Lyapunov equation 
\begin{equation}
B_{2}A_{2}+A_{2}^{T}B_{2}=-I,  \tag{5.12}
\end{equation}%
where 
\[
B_{2}=\left[ 
\begin{array}{ccc}
b_{11} & b_{12} & b_{13} \\ 
b_{21} & b_{22} & b_{23} \\ 
b_{31} & b_{32} & b_{33}%
\end{array}%
\right] . 
\]%
The equation $\left( 5.12\right) $ is reduced to the linear system of
algebraic equation in $b_{ij}$ and by solving this system we get

\[
b_{13}=0\text{, }b_{22}=-\frac{1}{2c_{22}},\text{ }b_{12}=b_{21}=-\frac{%
a_{21}\delta }{2c_{22}\left( c_{11}+c_{22}\right) },\text{ } 
\]%
\ 
\[
b_{11}=-\frac{1}{c_{11}}\left[ \frac{a_{21}^{2}\delta ^{2}}{2c_{22}\left(
c_{11}+c_{22}\right) }+\frac{1}{2}\right] ,\text{ }b_{23}=b_{32}=0,\text{ } 
\]%
\begin{equation}
b_{33}=-\frac{1}{2d_{33}},\text{ }b_{13}=b_{31}=0\text{, }B_{2}=\left[ 
\begin{array}{ccc}
b_{11} & b_{12} & 0 \\ 
b_{12} & b_{22} & 0 \\ 
0 & 0 & b_{33}%
\end{array}%
\right] .  \tag{5.q$_{2}$}
\end{equation}

Moreover, 
\begin{equation}
P_{2}\left( \lambda \right) =\left( b_{33}-\lambda \right) \left[ \lambda
^{2}-\left( b_{11}+b_{22}\right) \lambda -\left(
b_{12}^{2}-b_{11}b_{22}\right) \right] =0.  \tag{5.13}
\end{equation}

In view of the assumption (1) and (2) it is clear to see that%
\[
b_{11}>0,\text{ }b_{22}>0,\text{ }b_{33}>0,\text{ }b_{12}<0. 
\]%
By assumption (3),%
\[
\left( b_{11}+b_{22}\right) ^{2}+4\left( b_{12}^{2}-b_{11}b_{22}\right)
=\left( b_{11}-b_{22}\right) ^{2}+4b_{12}^{2}\geq 0\text{.} 
\]%
So, $\left( 5.13\right) $ has positive roots, i.e. the matrix $B_{2}$ is
positive defined for all $x$. Hence, the quadratic function 
\begin{equation}
V_{2}\left( x\right) =X^{T}P_{2}X=b_{11}x_{1}^{2}+2b_{12}x_{1}\left(
x_{2}-\delta \right) +2b_{22}\left( x_{2}-\delta \right) ^{2}+b_{33}x_{3}^{2}
\tag{$5.v_{2}$}
\end{equation}%
is a positive defined Lyapunov function candidate in certain neighborhood of 
$E_{2}\left( 0,\delta ,0\right) $. We need to determine a domain $\Omega
_{2} $ about the point $E_{2},$ where $\dot{V}_{2}\left( x\right) $ is
negative defined and a constant $C$ such that $\Omega _{C}$ is a subset of $%
\Omega _{2}$. By assuming $x\in R_{+}^{3}$ we will find the solution set of
the following inequality, 
\[
\dot{V}_{2}\left( x\right) =\dsum\limits_{k=1}^{3}\frac{\partial V_{2}}{%
\partial x_{k}}\frac{dx_{k}}{dt}= 
\]%
\[
\left[ 2b_{11}x_{1}+2b_{12}\left( x_{2}-\delta \right) \right] \left[ \left(
1-x_{1}\right) -a_{12}x_{2}-a_{13}x_{3}-g_{1}\left( \upsilon \right) \right]
x_{1}+ 
\]%
\[
\left[ 2b_{22}\left( x_{2}-\delta \right) +2b_{12}x_{1}\right] \left[
r_{2}\left( 1-x_{2}\right) -a_{21}x_{1}-g_{2}\left( \upsilon \right) \right]
x_{2}+ 
\]%
\begin{equation}
b_{33}x_{3}^{2}\left[ \frac{r_{3}x_{1}}{x_{1}+k_{3}}-a_{31}x_{1}-d_{3}-g_{3}%
\left( \upsilon \right) \right] <0.  \tag{5.14}
\end{equation}%
By assumptions, the coefficients of terms $x_{1}^{3},\ x_{2}^{3}$, $%
x_{1}x_{2}^{2},$ $x_{1}^{2}x_{2},$ $x_{1}x_{2}x_{3}$ and some coeffcients of 
$x_{1}x_{2}$, $x_{1}x_{3}$\ are negative. So, by assumption (2) and by (6q$%
_{2}$) for $x\in R_{+}^{3}$ the inequality $\left( 6.14\right) $ holds if

\[
\left[ 2b_{11}\left( 1-g_{1}\left( \upsilon \right) \right) -2\delta
b_{12}g_{1}\left( \upsilon \right) \right] x_{1}^{2}+2\delta b_{22}\left[
r_{2}+g_{2}\left( \upsilon \right) \right] x_{2}^{2}- 
\]%
\[
b_{33}\left( d_{3}+g_{3}\left( \upsilon \right) \right) x_{3}^{2}+\left[
-2\delta b_{12}a_{12}+2\delta b_{22}a_{21}-2b_{12}g_{1}\left( \upsilon
\right) \right] x_{1}x_{2}+ 
\]%
\[
2\delta b_{12}g_{1}\left( \upsilon \right) x_{1}+\left[ 2b_{22}g_{2}\left(
\upsilon \right) \left( \delta -r_{2}\right) \right] x_{2}+ 
\]%
\[
-\left[ 2b_{12}r_{2}+2b_{22}a_{21}\right] x_{1}x_{2}^{2}+\left\{ 2b_{12}%
\left[ g_{1}\left( \upsilon \right) -a_{21}\right] -2b_{11}a_{12}\right\}
x_{1}^{2}x_{2}<0. 
\]%
By assumption (3) 
\[
-\left[ 2b_{12}r_{2}+2b_{22}a_{21}\right] <0,\text{ }\left\{ 2b_{12}\left[
g_{1}\left( \upsilon \right) -a_{21}\right] -2b_{11}a_{12}\right\} <0. 
\]%
Hence, in view of inequalities $\left( 6.2\right) $, the \ above inequality
holds if 
\[
\mu _{0}\left[ x_{1}^{2}+\left( x_{2}-\delta \right) ^{2}+x_{3}^{2}\right]
<\mu _{0}\delta ^{2}-b_{12}\left( \delta a_{12}+g_{1}\left( \upsilon \right)
\right) +\delta b_{22}a_{21}, 
\]%
\begin{equation}
2\delta b_{12}g_{1}\left( \upsilon \right) x_{1}+\left[ 2b_{22}g_{2}\left(
\upsilon \right) \left( \delta -r_{2}\right) +2\mu _{0}\delta \right]
x_{2}\leq 0,  \tag{5.15}
\end{equation}%
where 
\[
\mu _{1}=\left[ 2b_{11}\left( 1-g_{1}\left( \upsilon \right) \right)
-2\delta b_{12}g_{1}\left( \upsilon \right) \right] ,\text{ }\mu
_{2}=2\delta b_{22}\left[ r_{2}+g_{2}\left( \upsilon \right) \right] , 
\]%
\[
\mu _{3}=b_{33}\left( d_{3}+g_{3}\left( \upsilon \right) \right) ,\text{ }%
\mu _{0}=\max \left\{ \mu _{1},\mu _{2},\mu _{3}\right\} .\text{ } 
\]

\bigskip From $\left( 5.q_{2}\right) $ it is easy to see that $\mu _{k}>0.$
Hence, $\dot{V}_{2}$ is negative defined on the domain $\Omega
_{2}=B_{r}\left( \bar{x}\right) \cap \Omega _{\delta },$ where%
\[
\bar{x}=\left( 0,\delta ,0\right) ,\text{ }r=\left( \frac{\eta }{\mu _{0}}%
\right) ^{\frac{1}{2}},\text{ }r\leq \sqrt{K_{1}^{2}+K_{2}^{2}+K_{3}^{2}},%
\text{ } 
\]%
\[
\eta =\mu _{0}\delta ^{2}-b_{12}\left( \delta a_{12}+g_{1}\left( \upsilon
\right) \right) +\delta b_{22}a_{21}+\mu _{2}\delta ^{2}, 
\]%
\begin{equation}
\Omega _{\delta }=\left\{ x\in R_{+}^{3}\text{, }2\delta b_{12}g_{1}\left(
\upsilon \right) x_{1}+\left[ 2b_{22}g_{2}\left( \upsilon \right) \left(
\delta -r_{2}\right) +2\mu _{0}\delta \right] x_{2}\leq 0\right\} . 
\tag{5.q$_{3}$}
\end{equation}%
i.e., the system $\left( 1.11\right) $ is global asymptotically\ stabile at $%
E_{1}\left( 0,\delta ,0\right) $.

\begin{center}
\bigskip \textbf{6. Basins of multiphase attraction sets }
\end{center}

\bigskip In this section we will derive momains of multipoint attraction
sets, $E_{0}\left( 0,0,0\right) ,$ $E_{1}\left( \gamma ,0,0\right) $, $%
E_{2}\left( 0,\delta ,0\right) ,$ where $\gamma $, $\delta $ were defined by 
$\left( 3.2\right) .$

We show in this section the following results

\textbf{Theorem 6.0. }Assume that all conditions of Theorem 5.0 are
satisfied. Then the basin of multiphase attraction set of the point $%
E_{0}\left( 0,0,0\right) $ belongs to the set $\Omega _{C}\subset \Omega
_{K},$ where the positive constant $C$ is defined in bellow.

\textbf{Proof. } Let us now find the set $\Omega _{C}\subset B_{r}\left( 
\bar{x}\right) ,$ where 
\[
C<\min_{\left\vert x\right\vert =r}V_{0}\left( x\right) =\lambda _{\min
}\left( P_{0}\right) r^{2}, 
\]%
here $P_{0}$ was defined by $\left( 6.0\right) $, $\lambda _{\min }\left(
P_{0}\right) $ denotes a minimum eignevalue of $P_{0}$, i.e. 
\[
\lambda _{\min }\left( P_{1}\right) =\min \left\{ b_{11},\text{ }b_{22},%
\text{ }b_{33}\right\} 
\]%
and 
\[
r\leq \sqrt{K_{1}^{2}+K_{2}^{2}+K_{3}^{2}.} 
\]

Moreover, for some $C>0$ the inclusion$\ \Omega _{C}\subset \Omega _{K}$
means the existence of $C>0$ so that $x\in \Omega _{C}$ implies $x\in \Omega
_{K}$, i.e. 
\[
0\leq x_{i}\leq K_{i}\text{, }K_{i}>\lambda _{\min }\left( P_{1}\right) . 
\]

\textbf{Theorem 6.1. }Assume that all conditions of Theorem 6.1 are
satisfied. Then the basin of multiphase attraction set of the point $%
E_{1}\left( \gamma ,0,0\right) $ belongs to the set $\Omega _{C}\subset
\Omega _{K}\cap \Omega _{\gamma }\cap B_{r}\left( \bar{x}\right) ,$ where $%
\Omega _{\gamma }$ was defined by $\left( 5.q_{1}\right) $, the positive
constant $C$ and $r$ were defined in bellow.

\textbf{Proof. }Let us now find the set $\Omega _{C}\subset B_{r}\left( \bar{%
x}\right) ,$ where 
\[
C<\min_{\left\vert x-\bar{x}\right\vert =r}V_{1}\left( x\right) =\lambda
_{\min }\left( P_{1}\right) \frac{\eta }{\mu _{0}},\text{ }\bar{x}=\left(
\gamma ,0,0\right) ,\text{ }P_{1}\left( \lambda \right) =\left\vert
B_{1}-\lambda \right\vert ,
\]%
here $B_{1}$ is a matrix defined by $\left( 5.5\right) $, $\lambda _{\min
}\left( P_{1}\right) $ denotes a minimum eignevalue of $P_{1}$, $\eta $, $%
\mu _{0}$\ were defined by $\left( 5.10\right) $ and $\left( 5.11\right) .$
Moreover, for some $C>0$ the inclusion$\ \Omega _{C}\subset \Omega _{\gamma }
$ means the existence of $C>0$ so that $x\in \Omega _{C}$ implies $x\in
\Omega _{\gamma }$, i.e. 
\[
0\leq x_{i}\leq K_{i}\text{,  }x_{1}\leq -\frac{1}{\alpha _{1}+2\mu _{1}}%
\left[ \alpha _{2}x_{2}+\alpha _{3}x\right] ,
\]%
$\alpha _{\imath },$ $\mu _{\imath }$ were defined by $\left( 5.q_{2}\right)
.$

So, 
\[
x\in B_{r}\left( \bar{x}\right) =\left\{ x\in R^{3}\text{, }\left\vert x-%
\bar{x}\right\vert <r_{0}\right\} ,
\]%
where 
\[
\text{ }r_{0}=\min \left\{ \left( \frac{\mu _{0}}{\mu _{1}\gamma ^{2}}%
\right) ^{\frac{1}{2}},\text{ }\left[ 2\gamma ^{2}+\left( 2\varkappa
_{2}^{2}+1\right) K_{2}^{2}+\varkappa _{2}^{2}K_{3}^{2}\right] ^{\frac{1}{2}%
}\right\} ,\text{  }
\]%
\[
\varkappa _{2}=\left\vert \frac{\alpha _{2}}{\alpha _{1}+2\mu _{1}}%
\right\vert ,\text{ }\varkappa _{3}=\left\vert \frac{\alpha _{3}}{\alpha
_{1}+2\mu _{1}}\right\vert .
\]

Then we obtain that 
\[
C<\min_{\left\vert x-\bar{x}\right\vert =\bar{r}}V_{1}\left( x\right)
=\lambda _{\min }\left( P_{1}\right) \bar{r}^{2},
\]%
i.e. 
\[
C<\lambda _{\min }\left( P_{1}\right) \bar{r}^{2},\text{ }r=\min \left\{
r_{0},\text{ }\bar{r}\right\} .
\]

Then we obtain that 
\[
C<\min_{\left\vert x-\bar{x}\right\vert =\bar{r}}V_{1}\left( x\right)
=\lambda _{\min }\left( P_{1}\right) \bar{r}^{2},
\]%
i.e. 
\[
C<\lambda _{\min }\left( P_{1}\right) \bar{r}^{2},\text{ }r=\min \left\{
r_{0},\text{ }\bar{r}\right\} .
\]

Now, we consider the equilibria point $E_{2}\left( 0,\delta ,0\right) $ and
prove the following result

\textbf{Theorem 6.2. }Assume all conditions of Theorem 5.2 are satisfied.
Then the basin of multiphase attraction set of $E_{2}\left( 0,\delta
,0\right) $ belongs to the set $\Omega _{C}\subset \Omega _{K}\cap \Omega
_{\delta }\cap B_{r}\left( \bar{x}\right) ,$ where $\Omega _{\delta }$ was
defined by $\left( 5.q_{3}\right) $ and 
\[
\text{ }\Omega _{C}=\left\{ x\in R^{3}\text{: }V_{2}\left( x\right) \leq C%
\text{ }\right\} ,\text{ }\bar{x}=\left( 0,\delta ,0\right) , 
\]%
here $V_{2}\left( x\right) $ was defined by $\left( 6.v_{2}\right) $,\ the
constants $C$ and $r$ are defined in bellow.

\textbf{Proof. }We will find $C>0$ such that $\Omega _{C}\subset \Omega
_{K}\cap B_{r}\left( \bar{x}\right) \cap \Omega _{\delta }$. It is clear to
see that $\Omega _{C}\subset B_{r}\left( \bar{x}\right) ,$ when 
\[
C<\min_{\left\vert x-\bar{x}\right\vert =r}V_{2}\left( x\right) =\lambda
_{\min }\left( P_{2}\right) r^{2},\text{ }\bar{x}=\left( 0,\delta ,0\right)
, 
\]%
here $\lambda _{\min }\left( P_{2}\right) $ denotes a minimum eigenvalue of $%
P_{2}$, i.e. 
\[
\lambda _{\min }\left( P_{2}\right) =\min \left\{ \frac{1}{2d_{3}},\text{ }%
\frac{\left( b_{22}+\frac{1}{2r_{2}}\right) \pm \sqrt{b_{11}^{2}+b_{12}^{2}+%
\frac{1}{4r_{2}^{2}}-\frac{b_{11}}{r_{2}}}}{2}\right\} . 
\]

Moreover, for some $C>0$ the inclusion$\ \Omega _{C}\subset \Omega _{\delta
} $ means the existence of $C>0$ so that $x\in \Omega _{C}$ implies $x\in
\Omega _{\delta }$, i.e.%
\begin{equation}
x\in \Omega _{K}\text{, }x_{2}\leq \beta x_{1},\text{ }\beta =-\frac{2\delta
b_{12}g_{1}\left( \upsilon \right) }{\left[ 2b_{22}g_{2}\left( \upsilon
\right) \left( \delta -r_{2}\right) +2\mu _{0}\delta \right] }.  \tag{6.1}
\end{equation}%
So, $x\in B_{\bar{r}}\left( \bar{x}\right) $, where $\mu _{i}$ were defined
by $\left( 5.15\right) $ and 
\[
\text{ }\bar{r}=\left[ K_{1}^{2}+\beta ^{2}\left( K_{1}-\delta \right)
^{2}+K_{3}^{2}\right] ^{\frac{1}{2}}. 
\]%
Then we obtain that 
\[
C<\min_{\left\vert x\right\vert =r_{0}}V_{2}\left( x\right) =\lambda _{\min
}\left( P_{2}\right) r_{0}^{2}, 
\]%
i.e. 
\[
C<\lambda _{\min }\left( P_{2}\right) r^{2}\text{ for }r=\min \left\{ r_{0},%
\text{ }\bar{r}\right\} . 
\]

\textbf{Remark 6.1. }It is clear to see that if $a_{21}\geq \frac{a_{12}-1}{%
r_{2}}$, then the assumption (3) is satisfied. Moreover, if $a_{12}+\frac{1}{%
2}\left( 1+a_{13}\right) >r_{2},$ then the assumption (4) holds.

\textbf{Remark 6.2. }The assumptions (3) and $\mu _{i}>0$ can be realized in
terms of coefficients of $\left( 1.3\right) $ by using $\left( 5.20\right) .$

\textbf{Conclusion. \ }1. Multipoint Initial Condition (MIC): The condition
determining initial state of complex system consisted of tumor cells, Immune
cells -- natural killer and host cells densities at the beginning of
observation; this condition is can be enriched with addition to the system
drug concentration. 

2. MIC is operated not only by changing the concentrations of the complex
system parameters, but also taking into account the impact of the drug
effects on the parameters, in order to direct the system to possible
equilibrium points. 

3. The possible equilibrium points can be selected and reached by the help
by changing the situation in behalf of a new, more effective equilibrium
point, one of multimodal attraction points, which multimodal attraction
basin is consisted of. This operation drifts the complex system into
situation, where Tumor Cells are trying to reach at least \textquotedblleft
Dormant State\textquotedblright\ that creates new chance to get more
important condition as healthy attraction point, which can be estimated as
globally stable condition. 

a) All these operations are traced beginning local lipschitz condition
(Theorem 2.1) application to IVP, being developed to global lipschitz
condition (Theorem 2.2). The conditions make it possible to figure out
positive multiphase attractor and multiphase basin; 

b) By the help of Remarks (3.1) and (3.2.) the equilibrium points are
described. These points relate to situation impeding tumor cell growth; 

c) Theorem (4.1) makes it possible to reach asymptotically stable point
(point, the function in the included conditions continuously approaches the
point) depending on relationships between kill rate of drug and rate of
tumor cell proliferation. Rate of Immune cell killing by tumor cells is also
considered for reaching the point. For reaching the target, condition of
positivity in relationships between combinations of different rates (cii)
were included;

d) By the help of Theorems (4.1. -- 4.7.), the main parametric conditions
are proved and determined for strong multipoint attractor, corresponding to
suppression of tumor growth, the possible healthy stable point in Multimodal
attraction basin; 

e) The important condition as a strong impact on tumor growth is reached by
proving theorems (5.0. -- 5.2) as global stability of equilibria points and
multiphase attraction sets;

(f) Theorems (6.0. -- 6.2.) describe the condition determining Multiphase
Attraction Set Basin, within which tumor cells are successfully suppressed;

Taking into account different and effective features of mathematical
modelling and its possibilities to figure out a problem in dynamics on the
basis of its logic properties, it was surely pointed out the characteristics
of a mathematical model to use in description of needed processes of a given
dynamic system with identified problems. In this paper, a three dimensional
model was devoted to mathematical description and regulation possibilities
of uncontrolled tumor processes by organism as a complex system. The
dynamics of interactions of the dimensions corresponded to tumor cells,
immune cells and healthy -- \textquotedblleft host\textquotedblright\ --
cells were given as forces of vectors, with addition another one vector as a
drug, negatively or positively converging to the basins of attractions,
depending on their importances for the complex system consisted of the four
factors. In order to make the model subjected to control, there was included
multiphase IVP, describing the system's important parameters to operate with
it in the farther processes of stages of development. The model was
undergone different changes to determine its limits of survival: it was
determined the conditions of boundedness the system can be restricted,
invariance in non- negativity, which means the model keeps its properties of
reactions to changing in proper way, being subjected to different analysis,
and the circumstances the system can be forced to be dissipated in. The
system was exposed to changing pressures to estimate its convenience to
biologically important properties as points of equilibria and local
stability conditions. The next step in exploring of the model were very
complex and logistic approaches to its properties for verification of the
conditions, providing the global equilibria points and multimodal attraction
sets, having biologically strong value in regulation of the processes
towards the positive effects of feasible medical external implementation at
the convenient stages, determined by multimodal attraction basins.

It is reasonable to observe of pationts, to get analysis (e.g. the status of
blood cells shown above) at different times that determines the states
patient at these times (maltipoint times) and its correlation. Moreover, to
define the system to be in stable, healthy state if it is in the basin of
multimodal attraction healthy point. Furthermore, it is reasonable that a
mathematical model of such a system should include at least two stable
multimodal attracting basin, one of which is considered \textquotedblleft
healthy\textquotedblright\ and another which is \textquotedblleft
diseased\textquotedblright . For if this were not the case, we would not
observe both types of behavior. A system with only an multipoint attracting
healthy state would never need to be treated, since it would naturally move
back to this state despite any exogenous shock. On the other hand, a system
with no multimodal attracting healthy state would never stay cured, and no
remission from disease would ever be observed.

\textbf{Acknowledgements}

The authors is thanking to Assoc. Prof. of the Department of Immunology of
Yeditepe University G\"{u}lderen Yanikkaya Demirel, Assoc. Prof. of the
Department of Medical Microbiology of Yeditepe University \.{I}brahim G.
Acuner and Prof. Dr. Faculty of the Helth Sciences of Okan University Aida
Sahmurova according to their valuable suggestions in the field of medicine
and biology.

\textbf{References}

\begin{enumerate}
\item Kuznetsov V. A., Makalkin I. A., Taylor M. A., Perelson A. S.,
Nonlinear dynamics of immunogenic tumors: parameter estimation and global
bifurcation analysis, Bull. Math. Biol. 1994(56), 295--321.

\item Adam J. A, Bellomo C., A survey of models for tumor-immune system
dynamics, Boston, MA: Birkhauser, 1996.

\item Eftimie R, Bramson J. L., Earn D. J. D., Interactions between the
immune system and cancer: a brief review of non-spatial mathematical models,
Bull. Math. Biol. 2011(73), 2--32.

\item Kirschner D, Panetta J., Modelling immunotherapy of the tumor--immune
interaction, J. Math. Biol. 1998(37), 235--52.

\item de Pillis L. G., Radunskaya A., The dynamics of an optimally
controlled tumor model: a case study, Math. Comput. Modell 2003(37),1221--44.

\item Nani F., Freedman H. A., Mathematical model of cancer treatment by
immunotherapy, Math Biosci. 2000(163), 159--99.

\item Altrock P. M., Lin L. Liu., Michor F., The mathematics of cancer:
integrating quantitative models, Nature Reviews Cancer. 15 (2015), No.12,
730-745.

\item Itik I. M., Banks S. P., Chaos in a three-dimensional cancer model,
Int J. Bifurcation Chaos 2010, 2010(20), 71--79.

\item Levine, H., A., Pamuk S., Sleeman B. D., Mathematical modeling of
capillary formation and development in tumor angiogenesis, Penetration into
the Stroma bulletin of Mathematical Biology (2001) 63, 801--863.

\item Starkov, K. E., Krishchenko A. P., On the global dynamics of one
cancer tumour growth model, Commun. Nonlinear Sci. Numer. Simul. 19 (2014),
1486--1495.

\item Wodarz D. and Komarova N. L., Dynamics of Cancer. Mathematical
Foundations of Oncology, World Scientific, Singapore (2014).

\item Jackson T., Komarova N. and Swanson K., Mathematical Oncology: Using
Mathematics to Enable Cancer Discoveries, American Mathematical Monthly,
121(9), (2014), 840-856.

\item Bellomo N., Preziosi L., Modelling and mathematical problems related
to tumor evolution and its interaction with the immune system, Math. Comput.
Modelling 32 (3/4), 413-452 (2000).

\item Khalil H. Nonlinear systems, NJ: Prentice Hall, 2002.

\item Verhulst, F., Nonlinear differential equations and dynamical systems,
Springer-Verlag Berlin Heidelberg 1996.

\item A. Swierniak and J. Smieja, Cancer chemotherapy optimization under
evolving drug resistance, Nonlinear Analysis- Theory, Methods and
Applications 47 (1)) 375-386 (2001).

\item D. Barbolosi and A. Iliadis, Optimizing drug regimens in cancer
chemotherapy: A simulation study using \symbol{126}1 pk-pd model, Computers
in Biology and Medicine 31 (3), 157-172 (2001).

\item A. Matveev and A. Savkin, Optimal chemotherapy regimens: Influence of
tumours on normal ceils and several toxicity constraints. IMA Journal of
Mathematics Applied in Medzcin,e and Bzology 18 (1). \%5-\symbol{126}40
(2001).

\item V. B. Shakhmurov, https://arxiv.org/submit/2126388, On the dynamics of
a cancer tumor growth model with multiphase structure, 2018
\end{enumerate}

\end{document}